\documentclass[smallextended]{svjour3m}       

\usepackage{multirow}

\newcommand{\grad}{\mathop{\rm grad}\nolimits}
\renewcommand{\div}{\mathop{\rm div}\nolimits}
\newcommand{\const}{\mathop{\rm const}\nolimits}

\smartqed  
\usepackage{graphicx}

\journalname{arXiv.Org}
\begin{document}

\title{Iterative methods for solving the pressure problem at multiphase filtration}

\titlerunning{Iterative methods for solving the pressure problem}       

\author{P. Vabishchevich and M. Vasil'eva}

\institute{P. Vabishchevich \at
              Nuclear Safety Institute, 
              52, B. Tulskaya, 115191 Moscow, Russia \\
              \email{vab@ibrae.ac.ru} \\
           M. Vasil'eva \at
              North-Eastern Federal University
              58, Belinskogo, 677000 Yakutsk, Russia \\
              \email{vasilyeva\_mv@mail.ru}
}

\date{Submitted to arXiv.org July 27, 2011}

\maketitle

\begin{abstract}
Applied problems of oil and gas recovery are studied numerically using the mathematical models of multiphase fluid flows in porous media. 
The basic model includes the continuity equations and the Darcy laws for each phase, as well as the algebraic expression for the sum of saturations.
Primary computational algorithms are implemented for such problems  using the pressure equation.
In this paper, we highlight the basic properties of the pressure problem and discuss the necessity of their fulfillment at the discrete level.
The resulting elliptic problem for the pressure equation is characterized by a non-selfadjoint operator.
Possibilities of approximate solving the elliptic problem are considered using the iterative methods.
Special attention is given to the numerical algorithms for calculating the pressure on parallel computers.

\end{abstract}

\section{Introduction}\label{s1}
Mathematical modeling of multicomponent flows in porous media is of
great importance in oil and gas recovery.
Traditionally, the hydrodynamic simulators for these applications are based on three-phase black oil model 
\cite{Aziz,Peaceman}.
A mathematical model of fluid dynamics in porous media includes differential equations, 
which express the conservation laws of mass and momentum \cite{Bear,Dagan}.
First of all, there are used the continuity equations describing the mass conservation law for
each separate phase.
The momentum  equations in a porous medium are written in the form of Darcy's law, 
which links the velocity with the pressure.
When capillary effects are omitted, the pressure is common to all phases.

Applied mathematical models of mass transfer processes in porous media 
are essentially nonlinear and difficult to study \cite{trangenstein1989mathematical,Vazquez}.
Next, it is necessary in these models to implement the closure of the system of equations 
via the constant sum of all saturations. 
Such algebraic components of the model should be taken into account
in constructing computational algorithms to predict multiphase flows in porous media. 
\cite{chen2006computational,fanchi2006principles}.

Two classes of methods are used to solve approximately unsteady boundary value problems 
for coupled systems of partial differential equations.
The first of them employs various implicit schemes for the initial system of equations.
In this case, we face some computational problems in the transition to a new time level.
The second class of methods reduces computational costs by means of using splitting schemes 
and solving simpler problems at the new time level --- 
splitting with respect to physical processes \cite{GIMarchuk_1990a,SamVabAdd}.
The above two classes of methods 
are presented through the fully implicit method (FIM) 
and the implicit pressure explicit saturation (IMPES) approach
\cite{Aziz,chen2006computational,Peaceman}.

FIM is widely used in hydrodynamic modeling of oil and gas recovery \cite{yeoh2010computational}.
The fully implicit approximation is used in FIM for all equations of the mathematical model.
It allows to expect stability of the method and possibility to use large time steps.
The basic drawback of the method is connected with its complexity -- we have to
solve a large system of nonlinear equations.

IMPES provides more efficient algorithms for solving the problem at each time level.
In this approach, we formulate the problem for the pressure with implicit approximations in time.
After evaluation of the pressure all other unknowns are calculated via explicit approximations.
Unfortunately, the problem of stability (time step restriction) is typical for the IMPES method.
Therefore, various modifications of IMPES have been developed in order to improve its stability.
For example, after evaluation of the pressure we can use
implicit approximations for the calculation of saturations \cite{watts1986compositional} (the sequential method).

In this paper, we highlight the main features of the pressure problem which should be 
taken into account in constructing computational algorithms.
There are discussed here possibilities of obtaining the pressure equation  --- 
elliptic for incompressible media and parabolic for compressible ones.
If the condition of the constant sum of saturations is treated explicitly then 
the corresponding elliptic operator of the pressure  problem is non-selfadjoint.
This fact should be taken into account in constructing iterative algorithms.

The paper is organized as follows.
A basic system of equations is formulated 
in Section \ref{s2} to describe multicomponent fluid flows in porous media. 
This mathematical  model is obtained at assumptions 
that capillary and gravity forces are negligible.
The pressure equation is derived in this section. 
The main features of the grid problem for the pressure are discussed in Section \ref{s3}.
The simplest uniform grids  for the problem in a rectangle are used.

The emphasis is on iterative methods for calculating the pressure at the new time level.
The two-dimensional test problem is described in Section \ref{s4}.
The possibility of using standard iterative methods with preconditioners is discussed 
in Section  \ref{s5}.
In section  \ref{s6}, we present the results of the iterative solving 
model  pressure problems on parallel computers.
Conclusions are summarized in Section \ref{s7}.

\section{The pressure problem}\label{s2}
In this section we formulate the basic mathematical model for fluid flows in porous media.
The system of governing equations for multicomponent flows includes  the continuity equation 
for each phase, where $\alpha =1,2,\dots, m$ --- the phase index.
The mass conservation law for each particular phase is expressed by the following equation
\begin{equation}\label{2.1}
  \frac{\partial (\phi \, b_\alpha S_\alpha)}{\partial t} +
  \div (b_\alpha  \mathbf{u}_\alpha) = 
  - b_\alpha q_\alpha,
  \quad \alpha =1,2,\dots, m .
\end{equation}
Here $\phi$ stands for the porosity,
$b_\alpha$ is the phase density,
$S_\alpha$ --- the phase saturation,
$\mathbf{u}_\alpha$ --- the velocity, and
$q_\alpha$  --- the volumetric mass source.

For simplicity, we neglect the capillary and gravitational forces.
In this simplest case the equation of fluid motion in porous media 
has the form of Darcy's law, where the velocity is directly determined by the common pressure:
\begin{equation}\label{2.2}
  \mathbf{u}_\alpha =
  - \frac{k_\alpha}{\mu_\alpha} \, \mathsf{k} \cdot \grad p ,
  \quad \alpha =1,2,\dots, m .
\end{equation}
In (\ref{2.2}),  $\mathsf{k}$ is the absolute permeability
(in general, symmetric second-rank tensor)
$k_\alpha$  --- the relative permeability,
$\mu_\alpha$ --- the phase viscosity and $p$ --- the pressure.

The unknowns in the system of equations (\ref{2.1}), (\ref{2.2})  
are the phase saturations $S_\alpha, \ \alpha =1,2,\dots, m$ 
and the pressure ($m+1$ unknowns in all).
In the simplest case, the coefficients in equations (\ref{2.1}), (\ref{2.2}) are
defined as some relations
\[
  \phi = \phi(p),
  \quad b_\alpha = b_\alpha(p),
  \quad q_\alpha = q_\alpha(S_\alpha),
  \quad k_\alpha = k_\alpha(S_\alpha),
  \quad \mu_\alpha = \const .
\]
For the sum of saturations of all phases we have
\begin{equation}\label{2.3}
  \sum_{\alpha =1}^{m} S_\alpha = 1.
\end{equation}
After substituting (\ref{2.2}) in (\ref{2.1}) and taking into account (\ref{2.3}), 
we have a system of $m+1$ equations for $m+1$ unknowns.

The system of equations (\ref{2.1})--(\ref{2.3}) is the basis for the description 
of multicomponent flows in porous media. In this system we have not any separate equation for the pressure.
Equations (\ref{2.1}) can be considered as the transport equation for each phase, 
whereas the algebraic relation (\ref{2.3}) can be treated as the equation for the pressure.

Let us consider more convenient forms of system (\ref{2.1})--(\ref{2.3}), 
which lead to the typical problems of mathematical physics for the pressure.
It should be noted that such equivalent formulations do exist only at the differential level.
At the discrete level such equivalence of formulations is broken even for linear problems.
So, the choice of the initial form of the equations is essential for calculations.

The most natural way to derive the equation for the pressure is the following.
Divide each equation (\ref{2.1}) by $\phi \, b_\alpha > 0$ and add them together, which gives
\begin{equation}\label{2.4}
  \left ( \sum_{\alpha =1}^{m} \frac{S_\alpha}{\phi \, b_\alpha}
  \frac{d (\phi \, b_\alpha) }{d p} \right ) 
  \frac{\partial p}{\partial t} =
  \sum_{\alpha =1}^{m} \frac{1}{\phi \, b_\alpha}
  \div  \left ( \frac{b_\alpha k_\alpha}{\mu_\alpha} \, \mathsf{k} \cdot \grad p \right ) -
  \frac{1}{\phi} \sum_{\alpha =1}^{m} q_\alpha .
\end{equation}
With the natural assumption for compressible fluids
\[
  \frac{d (\phi \, b_\alpha) }{d p}  > 0,
  \quad \alpha =1,2,\dots, m .
\]
equation (\ref{2.4}) for the pressure  is the standard parabolic equation of second order.
In particular, the maximum principle holds  for its solutions \cite{friedman1964partial}.

When using equation (\ref{2.4}), the basic system of equations for flows in porous media 
can include $m$ equations
\begin{equation}\label{2.5}
  \frac{\partial (\phi \, b_\alpha S_\alpha)}{\partial t} -
  \div  \left ( \frac{b_\alpha k_\alpha}{\mu_\alpha} \, \mathsf{k} \cdot \grad p \right ) =
  - b_\alpha q_\alpha,
\end{equation}
for $S_\alpha, \ \alpha =1,2,\dots, m$ and equation (\ref{2.4}) for $p$.
In this case equation (\ref{2.3}) is a consequence of (\ref{2.4}), (\ref{2.5}).
The second approach of common use 
is connected with employing relation (\ref{2.3}) instead of one of equations (\ref{2.5}).
For example,  equation (\ref{2.4}) is treated as the pressure equation, 
equations (\ref{2.5}) are used for $S_\alpha, \ \alpha =1,2,\dots, m-1$  
whereas from  (\ref{2.3}) we get $S_m$
\begin{equation}\label{2.6}
  S_m = 1 - \sum_{\alpha =1}^{m-1} S_\alpha.
\end{equation}
Note that the above forms of equations for multicomponent flows in porous media 
are algebraically equivalent only at the differential level.
We will try to preserve the main points of this equivalence
at the discrete level \cite{vab2010}.

In the case of variable coefficients  $\phi \, b_\alpha$, the elliptic operator 
for the pressure equation (\ref{2.4}) is non-selfadjoint.
This fact leads to some problems in using implicit schemes 
for equation (\ref{2.4}). That is why some modifications are employed for the pressure equation.
For instance, we can obtain the pressure equation via the direct summation of equations (\ref{2.1})
taking into account equation (\ref{2.2})
\begin{equation}\label{2.7}
  \sum_{\alpha =1}^{m} 
  \frac{\partial (\phi \, b_\alpha S_\alpha)}{\partial t} =
  \sum_{\alpha =1}^{m} 
  \div  \left ( \frac{b_\alpha k_\alpha}{\mu_\alpha} \, \mathsf{k} \cdot \grad p \right ) -
  \sum_{\alpha =1}^{m} b_\alpha q_\alpha .
\end{equation} 
In this case we have the selfadjoint elliptic operator for the pressure.
However, this approach has some drawbacks.
In particular, this system of equations (\ref{2.1}), (\ref{2.2}),(\ref{2.7}) is not closed, 
because the basic algebraic relation (\ref{2.3}) is not involved.
Also, we can not say anything about the parabolic property of the pressure operator in equation (\ref{2.7}).

\section{The properties of the grid operators}\label{s3}
Let us consider  stationary and unsteady model problems, 
which are linear prototypes for the pressure problem in modeling multiphase flows.
Consider the two-dimensional problem in the rectangle
\[
  \Omega = \{ \ \mathbf{x} \ | \ \mathbf{x} = (x_1, x_2), 
  \ 0 < x_{\beta} < l_{\beta}, \ \beta =1,2 \}.
\]

In accordance with (\ref{2.4})  we solve in $\Omega$ the boundary problem for the equation
\begin{equation}\label{3.1}
   \frac{\partial u}{\partial t} +
   \sum_{\alpha =1}^{m} a_\alpha({\bf x}) \mathcal{L}_\alpha u = f({\bf x},t) ,
   \quad {\bf x}\in \Omega,
   \quad 0 < t < T,
\end{equation}
where $a_\alpha({\bf x}) \geq \varrho_\alpha, \ \varrho_\alpha  > 0, \ \alpha = 1,2, ..., m$, 
and elliptic operators $\mathcal{L}_\alpha$ are defined by
\begin{equation}\label{3.2}
   \mathcal{L}_\alpha u = - \frac{\partial }{\partial x_1} 
   \left (k_\alpha ({\bf x}) \frac{\partial u}{\partial x_1} \right ) 
   - \frac{\partial }{\partial x_2} \left (k_\alpha ({\bf x}) \frac{\partial u}{\partial x_2} \right ) ,
   \quad  \alpha = 1,2, ..., m,
\end{equation} 
under the standard assumptions $0 < \kappa_\alpha \leq k_\alpha  \leq \overline{\kappa}_\alpha$.
This equation is supplemented with homogeneous Dirichlet boundary conditions
\begin{equation}\label{3.3}
   u({\bf x},t) = 0,
   \quad {\bf x}\in \partial \Omega,
   \quad t > 0.
\end{equation}
In addition, the initial condition is given in the following form
\begin{equation}\label{3.4}
   u({\bf x},0) = u^0({\bf x}),
   \quad {\bf x}\in \Omega.
\end{equation}

In some cases (incompressible media) it is reasonable to consider the stationary problem.
The boundary value problem is formulated for the equation
\begin{equation}\label{3.5}
   \mathcal{A} u = f({\bf x}),
   \quad \mathcal{A} = 
   \sum_{\alpha =1}^{m} a_\alpha({\bf x}) \mathcal{L}_\alpha  ,
   \quad {\bf x}\in \Omega,
\end{equation}
which is supplemented by the boundary conditions (\ref{3.3}).

The approximate solution is given at the nodes of the uniform rectangular grid in  $\Omega$:
\[
   \bar{\omega} = \{ \mathbf{x} \ | \ \mathbf{x} = (x_1, x_2),
   \ x_\beta = i_\beta h_\beta,
   \ i_\beta = 0,1,...,N_\beta,
   \ N_\beta h_\beta = l_\beta, 
   \ \beta=1,2\} 
\]
and let $\omega$ be the set of internal nodes ($\bar{\omega} = \omega \cup \partial \omega$). 
For grid functions $y(\mathbf{x}) = 0, \ \mathbf{x} \in \partial \omega$ 
we define the Hilbert space $H = L_2({\omega})$  with the inner product and norm
\[
  (y,w) \equiv \sum_{{\bf x} \in \omega}
  y({\bf x}) w({\bf x}) h_1 h_2,
  \quad \|y\| \equiv (y,y)^{1/2} .
\]

Approximation in space for  problem (\ref{3.1})--(\ref{3.4}) will be performed
at the assumption that the coefficients and solution are sufficient smooth.
For the elliptic operator $\mathcal{L}_\alpha$ we put into the correspondence the grid operator
$\Lambda_{\alpha}$:
\[
  \Lambda_{\alpha} y =
  - \frac{1}{h_1^2} k_{\alpha}(x_1+0.5h_1,x_2)
  (y(x_1+h_1,x_2) - y(x_1,x_2))
\]
\[
  + \frac{1}{h_1^2} k_{\alpha}(x_1-0.5h_1,x_2)
  (y(x_1,x_2) - y(x_1-h_1,x_2)) 
\]
\[
  - \frac{1}{h_2^2} k_{\alpha}(x_1,x_2+0.5h_2)
  (y(x_1,x_2+h_2) - y(x_1,x_2))
\]
\begin{equation}\label{3.6}
  + \frac{1}{h_2^2} k_{\alpha}(x_1,x_2-0.5h_2)
  (y(x_1,x_2) - y(x_1,x_2-h_2)) ,
  \quad {\bf x}\in \omega 
\end{equation}
for all $\alpha =1,2,...,m$.
In $H$ the operators $\Lambda_{\alpha}, \ \alpha =1,2,...,m$ are selfadjoint 
and positive definite \cite{2001Samarskii,SamNik}:
\begin{equation}\label{3.7}
  \Lambda_{\alpha} = \Lambda_{\alpha}^* ,
  \quad \kappa_{\alpha} \delta E \leq  \Lambda_{\alpha} \leq \Delta \overline{\kappa}_{\alpha} E,
  \quad  \alpha = 1,2, ..., m,
\end{equation}
where $ E $ --- the identity operator, and
\[
  \delta = \delta_1+\delta_2,
  \quad \delta_{\beta} = 
  \frac{4}{h^2_{\beta }} \sin^2 \frac{\pi h_{\beta}}{2 l_{\beta}} ,
\]
\[
  \Delta = \Delta_1+\Delta_2,
  \quad \Delta_{\beta} = 
  \frac{4}{h^2_{\beta}} \cos^2 \frac{\pi h_{\beta}}{2 l_{\beta}} ,
  \quad \beta = 1,2 .
\]

The grid operator for the pressure problem can be represented as
\begin{equation}\label{3.8}
  A = \sum_{\alpha =1}^{m} a_\alpha({\bf x}) \Lambda_\alpha ,
  \quad {\bf x}\in \omega .
\end{equation} 
In general (non-constant coefficients $a_\alpha({\bf x}), \ \alpha = 1,2, ..., m$) 
the operator $A$ is non-selfadjoint.
It approximates the corresponding differential operator with the error of 
$\mathcal{O} (|h|^2)$, where $|h|^2 = h_1^2 + h_2^2$.

After discretization in space we go from (\ref{3.1})--(\ref{3.4}) 
to the differential-operator equation
\begin{equation}\label{3.9}
  \frac {dy}{dt} + A y = f(t),
  \quad 0 < t < T,
\end{equation}
considered on the set of grid functions $y(t) \in H$.
The initial condition is taken in the form
\begin{equation}\label{3.10}
  y(0) = u^0.
\end{equation}
For the stationary problem (\ref{3.3}), (\ref{3.5}) the grid analog has the form
\begin{equation}\label{3.11}
  A y = f .
\end{equation}

To solve approximately problem (\ref{3.9}), (\ref{3.10}), we use the standard two-level schemes.
Let $\tau$ be the fixed time step and
$y^n = y(t^n), \ t^n = n \tau$, $n = 0,1, ..., N, \ N\tau = T$.
Equation (\ref{3.9}) is approximated by the two-level scheme with weights
\begin{equation}\label{3.12}
  \frac{y^{n+1} - y^{n}}{\tau }
  + A(\sigma y^{n+1} + (1-\sigma) y^{n}) = \varphi^n,
  \quad n = 0,1, ..., N-1,
\end{equation}
where, for example, $\varphi^n = f(\sigma t^{n+1} + (1-\sigma) t^{n})$.
It is supplemented by the initial condition
\begin{equation}\label{3.13}
  y^0 = u^0 .
\end{equation}
The difference scheme  (\ref{3.12}), (\ref{3.13}) has the approximation error in time
$\mathcal{O} (\tau^2 + (\sigma - 0.5) \tau)$.

If we employ the fully implicit scheme ($\sigma = 1$), then
the transition to the new time level is performed through solving the grid problem
\begin{equation}\label{3.14}
  \left ( \frac{1}{\tau } E +  A \right ) y = f .
\end{equation}
The main subject of our consideration is the methods of solving grid problems  
(\ref{3.11}) and (\ref{3.14}), which are linear prototypes for stationary and unsteady problem 
for the pressure. The primary question here is the non-selfadjoint property of the grid operator $A$.

For the grid problem (\ref{3.11}) the maximum principle holds  \cite{2001Samarskii}.
With regard to considered approximations on the five-point stencil, 
we formulate it as follows \cite{SamVabConv}.
Consider the difference equation
\[
  \gamma({\bf x}) y({\bf x}) - \alpha_1({\bf x}) y(x_1-h_1,x_2) -
  \beta_1 ({\bf x}) y(x_1+h_1,x_2) -
\]
\begin{equation}\label{3.15}
  - \alpha_2({\bf x})  y(x_1,x_2-h_2) -
  \beta_2 ({\bf x}) y(x_1,x_2+h_2) =
  \varphi({\bf x}),
  \quad {\bf x} \in \omega ,
\end{equation}
which is supplemented by boundary conditions
\begin{equation}\label{3.16}
   y({\bf x}) = 0,
  \quad {\bf x} \in \partial\omega .
\end{equation}
We assume, that the coefficients of the difference scheme (\ref{3.15}) satisfy the conditions
\begin{equation}\label{3.17}
  \alpha_j ({\bf x}) > 0,
  \quad \beta_j ({\bf x}) > 0,
  \quad j = 1,2,
  \quad \gamma ({\bf x}) > 0,
  \quad {\bf x} \in \omega .
\end{equation}

Let in the difference scheme (\ref{3.15})-(\ref{3.16}) we have
$\varphi ({\bf x}) \ge 0$ for all ${\bf x} \in \omega$  
(or $\varphi ({\bf x}) \le 0$ for  ${\bf x} \in \omega$).
Then for
\begin{equation}\label{3.18}
  \gamma ({\bf x}) \ge
  \alpha_1 ({\bf x}) + \alpha_2 ({\bf x}) +
  \beta_1 ({\bf x}) + \beta_2 ({\bf x}),
  \quad {\bf x} \in \omega
\end{equation}
we have (the grid maximum principle)
$y ({\bf x}) \ge 0$ for all ${\bf x} \in \omega$ ($y({\bf x}) \le 0$ for ${\bf x} \in \omega$).
In our case (see (\ref{3.6}), (\ref{3.8})) fulfillment of the sufficient conditions 
(\ref{3.18}) can be verified directly.
Because of this, for the grid operator at the new time level (\ref{3.14}) we have the strict diagonal dominance.

To study properties of operators  $\mathcal{A}$ and $A$
in Hilbert spaces $\mathcal{H} = L_2(\Omega)$ and $H = L_2(\omega)$, it is convenient 
to treat $\mathcal{A}$ and $A$ as the corresponding convection-diffusion operators.
In this case it is possible to employ in our research the results from \cite{MortonCon,SamVabConv}.

Taking into consideration (\ref{3.2}) and (\ref{3.5}), we have the representation
\begin{equation}\label{3.19}
  \mathcal{A} = \sum_{\alpha =1}^{m} \mathcal{A}_\alpha,
  \quad \mathcal{A}_\alpha = \mathcal{D}_\alpha  + \mathcal{C}_\alpha, 
  \quad \alpha = 1,2, ..., m, 
\end{equation} 
where
\begin{equation}\label{3.20}
  \mathcal{D}_\alpha u  = - \div (d_\alpha ({\bf x}) \grad u),
\end{equation} 
\begin{equation}\label{3.21}
  \mathcal{C}_\alpha u  = \textbf{w}_{\alpha} \grad u .
\end{equation}   
The effective diffusion coefficient and convection velocity for the separate phase $\alpha$ are
\[
  d_\alpha =  k_\alpha a_\alpha,
  \quad \textbf{w}_{\alpha}= k_{\alpha} \grad  a_\alpha .   
\] 
Then the pressure operator takes the form of convection-diffusion operator 
with the convective term in the non-divergent form.
Note that application of equation (\ref {2.7}) to evaluate the pressure
corresponds to using only the diffusion part (\ref {3.20}) of  the operator (\ref{3.19}).

Operators of diffusion in the above assumptions about the coefficients
are self-adjoint and positive definite
in $\mathcal{H} = L_2(\Omega)$.
Next, we present some facts about the properties of convective transport operators.
A detailed discussion of these issues is given in the book \cite{SamVabConv}.

We have the following representation
\begin{equation}\label{3.22}
  \mathcal{C}_\alpha  = \overline{\mathcal{C}}_\alpha  -
  \frac{1}{2} \div \textbf{w}_{\alpha} ,
  \quad \alpha = 1,2, ..., m,   
\end{equation} 
where $\mathcal{C}_\alpha$ is the operator of convective transport in the symmetric form:
\begin{equation}\label{3.23}
  \overline{\mathcal{C}}_\alpha u =
  \frac{1}{2} \left ( \textbf{w}_{\alpha} \grad u + \div(\textbf{w}_{\alpha} u) \right ). 
\end{equation}
The operator $\overline{\mathcal{C}}_\alpha$ is skew-symmetric  in $\mathcal{H}$:
\begin{equation}\label{3.24}
  \overline{\mathcal{C}}_\alpha = - \overline{\mathcal{C}}_\alpha^* 
\end{equation} 
for any $\textbf{w}_{\alpha}({\bf x}), \ {\bf x} \in \Omega$.

From (\ref{3.22}) and (\ref{3.24}) we directly obtain the estimate for the energy 
of the convective transport operator $\mathcal{C}_\alpha$:
\begin{equation}\label{3.25}
  |(\mathcal{C}_\alpha u, u) | \leq \mathcal{M}_{\alpha} \| u\|^2 ,
\end{equation} 
\begin{equation}\label{3.26}
   \mathcal{M}_{\alpha} = \frac{1}{2} \| \div \textbf{w}_{\alpha} \|_{C(\Omega)},
   \quad \| u \|_{C(\Omega)} \equiv  \max_{{\bf x} \in \Omega} | u({\bf x})| .
\end{equation} 
It is interesting to consider the subordination estimate for the operator of convective
transport with respect to the diffusion operator.
In our model two-dimensional problem the corresponding estimate has the following form
\begin{equation}\label{3.27}
  \| \mathcal{C}_\alpha u\|^2 \leq \overline{\mathcal{M}}_{\alpha}  (\mathcal{D}_\alpha u, u) ,    
\end{equation} 
at $\textbf{w}_{\alpha} = (w_\alpha^{(1)}, w_\alpha^{(2)})$ with constant
\begin{equation}\label{3.28}
  \overline{\mathcal{M}}_{\alpha}  \leq \frac{2}{\varrho_\alpha  \kappa_\alpha}   
  \max_{\beta=1,2} \left \{ \left \|  
  \left ( w_\alpha^{(\beta)} \right )^2  \right \|_{C(\Omega)} \right \} .       
\end{equation}

These properties of differential operators of diffusion and convection 
(\ref{3.22}), (\ref{3.24}), (\ref{3.25}), (\ref{3.27}) are inherited not only for
the difference operators on rectangular grids \cite{SamVabConv}, 
but also for difference operators on irregular grids with the Delaunay triangulation \cite{vabVAGO}.
We consider this issue here for the grid operator (\ref{3.6}), (\ref{3.8}).
First of all, we are interested in the grid analog of  (\ref{3.19})--(\ref{3.21}).

Taking into account
\[
  - a(x) \frac{1}{h} \left ( k(x+h/2) \frac{y(x+h) - y(x)}{h} -
  k(x-h/2) \frac{y(x) - y(x-h)}{h} \right )  =
\]
\[
  - \frac{1}{h} \left ( \frac{a(x+h) + a(x)}{2} k(x+h/2) \frac{y(x+h) - y(x)}{h} - \right .
\]
\[
  \left .  
  \frac{a(x) + a(x-h)}{2} k(x-h/2) \frac{y(x) - y(x-h)}{h} \right )  +
\]
\[
  \frac{1}{2} \frac{a(x+h) - a(x)}{h} k(x+h/2) \frac{y(x+h) - y(x)}{h} +
\]
\[
  \frac{1}{2} \frac{a(x) - a(x-h)}{h} k(x-h/2) \frac{y(x) - y(x-h)}{h} ,
\]
similarly (\ref{3.19}) we obtain
\begin{equation}\label{3.29}
  A = \sum_{\alpha =1}^{m} A_\alpha,
  \quad A_\alpha = D_\alpha  + C_\alpha, 
  \quad \alpha = 1,2, ..., m .
\end{equation} 
The grid diffusion operator has the form
\[
  D_{\alpha} y =
\]
\[
  - \frac{a_{\alpha}(x_1+h_1,x_2) + a_{\alpha}(x_1,x_2)}{2h_1^2} k_{\alpha}(x_1+0.5h_1,x_2)
  (y(x_1+h_1,x_2) - y(x_1,x_2))
\]
\[
  + \frac{a_{\alpha}(x_1,x_2)+a_{\alpha}(x_1-h_1,x_2)}{2h_1^2} k_{\alpha}(x_1-0.5h_1,x_2)
  (y(x_1,x_2) - y(x_1-h_1,x_2)) 
\]
\[
  - \frac{a_{\alpha}(x_1,x_2+h_2) + a_{\alpha}(x_1,x_2)}{2h_2^2} k_{\alpha}(x_1,x_2+0.5h_2)
  (y(x_1,x_2+h_2) - y(x_1,x_2))
\]
\begin{equation}\label{3.30}
  + \frac{a_{\alpha}(x_1,x_2)+a_{\alpha}(x_1,x_2-h_2)}{2h_2^2} k_{\alpha}(x_1,x_2-0.5h_2)
  (y(x_1,x_2) - y(x_1,x_2-h_2)) .
\end{equation}
Similarly (\ref{3.7}) in $H = L_2(\omega)$ we have
\begin{equation}\label{3.31}
  D_{\alpha} = D_{\alpha}^* ,
  \quad  D_{\alpha} \geq \rho_\alpha \kappa_{\alpha} \delta E .
\end{equation}

Approximation of the convective part of the grid operator $A$ is conducted via setting 
the coefficients (effective velocity $\textbf{w}_{\alpha}$) on the grids shifted 
in the corresponding direction on the half-step.
Let us define with accuracy of $\mathcal{O}(|h|^2)$ the components of 
the grid analog of $\textbf{w}_{\alpha}$ using the following relations
\[
  w_\alpha^{(1)} (x_1+0.5h_1, x_2) = 
  \frac{a_{\alpha}(x_1+h_1,x_2) - a_{\alpha}(x_1,x_2)}{h_1} k_{\alpha}(x_1+0.5h_1,x_2) ,
\]
\begin{equation}\label{3.32}
  w_\alpha^{(2)} (x_1, x_2+0.5h_2) = 
  \frac{a_{\alpha}(x_1,x_2+h_2) - a_{\alpha}(x_1,x_2)}{h_2} k_{\alpha}(x_1,x_2+0.5h_2) . 
\end{equation}  
The convective transport operator in  representation (\ref{3.29}) has the form
\[
  C_\alpha y = \frac{1}{2} w_\alpha^{(1)} (x_1+0.5h_1, x_2) \frac{y(x_1+h_1,x_2) - y(x_1,x_2)}{h_1}  +
\]
\[
  \frac{1}{2} w_\alpha^{(1)} (x_1-0.5h_1, x_2) \frac{y(x_1,x_2) - y(x_1-h_1,x_2)}{h_1}  +
\]
\[
  \frac{1}{2} w_\alpha^{(2)} (x_1, x_2+0.5h_2) \frac{y(x_1,x_2+h_2) - y(x_1,x_2)}{h_2}  +
\]
\begin{equation}\label{3.33}
  \frac{1}{2} w_\alpha^{(2)} (x_1, x_2-0.5h_2) \frac{y(x_1,x_2) - y(x_1,x_2-h_2)}{h_2}  .
\end{equation}

The grid analogue of (\ref{3.11}) can be written as
\begin{equation}\label{3.34}
  C_\alpha = \overline{C}_\alpha -   \frac{1}{2} \div_h \textbf{w}_{\alpha} ,
\end{equation} 
where
\[
  \div_h \textbf{w}_{\alpha} = 
  \frac{w_\alpha^{(1)} (x_1+0.5h_1, x_2) - w_\alpha^{(1)} (x_1-0.5h_1, x_2)}{h_1} +
\]
\begin{equation}\label{3.35}
  \frac{w_\alpha^{(2)} (x_1+, x_20.5h_2) - w_\alpha^{(2)} (x_1, x_2-0.5h_2)}{h_2} . 
\end{equation} 
For the skew-symmetric part
\begin{equation}\label{3.36}
  \overline{C}_\alpha = - \overline{C}_\alpha^*
\end{equation} 
we have
\[
 \overline{C}_\alpha y = \frac{1}{2 h_1} w_\alpha^{(1)} (x_1+0.5h_1, x_2) y(x_1+h_1,x_2) -
\]
\[
 \frac{1}{2 h_1} w_\alpha^{(1)} (x_1-0.5h_1, x_2) y(x_1-h_1,x_2) +
\]
\[
 \frac{1}{2 h_2} w_\alpha^{(2)} (x_1, x_2+0.5h_2) y(x_1,x_2+h_2) -
\]
\begin{equation}\label{3.37}
 \frac{1}{2 h_2} w_\alpha^{(2)} (x_1, x_2-0.5h_2) y(x_1,x_2-h_2) .
\end{equation} 
The following grid analog of (\ref{3.25}) takes place:

\begin{equation}\label{3.38}
  |(C_\alpha y, y) | \leq M_{\alpha} \| y\|^2 ,
\end{equation} 
where now (see (\ref{3.26}))
\begin{equation}\label{3.39}
   M_{\alpha} = \frac{1}{2} \| \div_h \textbf{w}_{\alpha} \|_{C(\omega)},
   \quad \| y \|_{C(\omega)} \equiv  \max_{{\bf x} \in \omega} | y({\bf x})| .
\end{equation} 
The subordination inequality (see (\ref{3.27}) and (\ref{3.28})) has the form
\begin{equation}\label{3.40}
  \|(C_\alpha y\|^2 \leq \overline{M}_{\alpha} (D_\alpha y, y) ,
\end{equation} 
where
\begin{equation}\label{3.41}
  \overline{M}_{\alpha} = 
  \max\left \{ \left \| \left ( w_\alpha^{(1)} (x_1\pm0.5h_1, x_2)\right )^2  \right \|_{C(\omega)},
  \left \| \left ( w_\alpha^{(2)} (x_1, x_2\pm0.5h_2)\right )^2  
  \right \|_{C(\omega)} \right \} .       
\end{equation} 
The fundamental issue here is that for these approximations the constants $M_\alpha$ 
and $\overline{M}_\alpha$ are the complete grid analogues of the corresponding 
constants $\mathcal{M}_\alpha$ and $\overline{\mathcal{M}}_\alpha$ for the differential problem.

\section{The test problem}\label{s4}
Capabilities of iterative methods for approximate solving the pressure equation 
in modeling multiphase flows in porous media are illustrated here using the test grid problem.
We consider equation (\ref{3.14}), which corresponds to the calculation 
of one time step in the numerical solution of problem (\ref{3.1})--(\ref{3.4}).
Numerical experiments are conducted for problem (\ref{3.14}) with $f = 1$ in the 
unit square ($l_\beta = 1, \ \beta = 1,2$) on the grid $h = h_1 = h_2$ ($N = N_1 = N_2$).

Particular attention should be given to the coefficients of equations 
(\ref{3.1}), (\ref{3.2}) in order to take into account peculiarities of these problems, namely,
inhomogeneity of
$a_\alpha(\mathbf{x}), \ \alpha=1,2, ..., m$.
Taking into account (\ref{2.4}), we set
\[
   k_\alpha(\mathbf{x}) \sim \frac{1}{a_\alpha(\mathbf{x})} ,
   \quad  \alpha=1,2, ..., m .
\]
Consider two-phase medium ($m=2$) with an incompressible fluid as the first phase
\[
  a_1(\mathbf{x}) = 1,
  \quad k_1(\mathbf{x}) = 1 .
\]
Compressibility of the second phase is defined as follows:
\[
  a_2(\mathbf{x}) = \exp(-\xi ((x_1 - 0.5)^2 + (x_2 - 0.5)^2)),
\]
\[
  k_2(\mathbf{x}) = \eta \exp(\xi ((x_1 - 0.5)^2 + (x_2 - 0.5)^2)) .
\]
The diffusion part of operator (\ref{3.20}) is 
\[
  \mathcal{D}_1 u  = - \div \grad u,
  \quad \mathcal{D}_2 u  = - \eta \div \grad u .
\]

Properties of the considered problems are defined(see (\ref{3.21})) 
by the vectors $\mathbf{w}_\alpha, \ \alpha=1,2, ...,m$.
For the test problem we have
\[
 \mathbf{w}_1 = 0,
 \quad \mathbf{w}_2 = (-2 \eta \xi (x_1 - 0.5), \ -2 \eta \xi (x_2 - 0.5) ) .
\]
In this case
\[
  \div \mathbf{w}_2 = - 4 \eta \xi .
\]
For the constants in the estimates  (\ref{3.25}) and (\ref{3.27}) we obtain
\[
   \mathcal{M}_2 =  2 \eta | \xi |,
   \quad \overline{\mathcal{M}}_2 =  2 \eta \xi^2 .
\]
Thus, the governing numerical parameters for this problem are $\eta$ and $\xi$.
The sign of $\xi$ can be any, moreover, it defines the fundamental difference in
the behavior of the solution (the pressure) in the vicinity of the production or injection well.

\section{Iterative solution of the problem}\label{s5} 
For numerical solving the test problem we use iterative methods.
In the corresponding grid equation (\ref{3.14}) the operator $A$ is non-selfadjoint.
Therefore, we use iterative methods for grid problems with unsymmetric matrices 
\cite{saad2003iterative,SamNik}.
The standard Generalized Minimal Residual Method 
(GMRES) with different preconditioners has been employed.

To solve the test problem, the  PETSc library \cite{PETSc} has been used. 
The PETSc library, developed in the Argonne National Laboratory, is a powerful set of 
freely available multi-platform compatible tools for the solution of 
large-scale problems governed by partial differential equations. 
Experiments were carried out with the following preconditioners:
\begin{description}
 \item[none] --- without preconditioning;
 \item[jacobi] --- the Jacobi method;
 \item[sor] --- the successive overrelaxation method;
 \item[ilu] --- the incomplete LU factorization;
 \item[mg] ---  the multigrid method.
 \end{description}

Table \ref{tab:xieta} shows the dependence of the computational cost
(the number of iterations) on the physical parameters of the problem.
Features of the problem are clearly defined by the parameters  $\eta$ and $\xi$.
The calculations were performed using the unpreconditioned GMRES method 
on the grid with $256 \times 256$ unknowns.
We see that with increasing of $\eta$ and/or $|\xi|$ the number of iterations decreases.
The same is true for negative values of $\eta$.
\begin{table}[!h]
\begin{center}
\caption{The dependence of the number of iterations on the physical parameters ($\xi$, $\eta$)}
\begin{tabular}{cccccc}
\hline
\multirow{2}{*}{$\eta$}  & \multicolumn{5}{c}{$\xi$}\\ 
              &$-10$ & $-1$ & $0.0$ & $1$  & $10$ \\
\hline $0.01$ & 4710 & 4686 & 4682  & 4678 & 4629 \\
$0.1$  & 4413 & 4726 & 4699  & 4652 & 3828 \\
$1$    & 1445 & 4790 & 4790  & 4289 & 1689 \\
$10$   & 879  & 4568 & 4884  & 3919 & 1088 \\
$100$  & 857  & 4493 & 4903  & 3856 & 1026 \\
\hline
\end{tabular} 
\label{tab:xieta}
\end{center}
\end{table}

Effect of  preconditioning on different grids is shown in Table \ref{tab:pc}. 
Calculations were performed at $\eta = 1$.
It is easy to see that the multigrid preconditioner is the best.
\begin{table}[!h]
\begin{center}
\caption{The number of iterations for different preconditioners depending on the grid}
\begin{tabular}{ccccccc}
\hline 
\multirow{2}{*}{grid} & \multirow{2}{*}{preconditioner} & \multicolumn{5}{c}{$\xi$} \\ 
 &    &$-10$ & $-1$ &$0.0$ & $1$  & $10$ 
\\ \hline
\multirow{5}{*}{$128 \times 128$}
&none   				& 538  & 1283 & 1270 & 1150 & 523  \\ 
&jacobi 				& 507  & 1283 & 1270 & 1150 & 517  \\ 
&sor    				& 222  & 280  & 284  & 281  & 156  \\ 
&ilu    				& 175  & 217  & 214  & 215  & 128  \\ 
&mg     				& 5    & 5    & 5    & 5    & 5    \\  
\hline 
\multirow{5}{*}{$256 \times 256$} 						   
&none   				& 1445 & 4790 & 4790 & 4289 & 1689 \\ 
&jacobi 				& 1443 & 4789 & 4790 & 4284 & 1675 \\ 
&sor    				& 350  & 807  & 765  & 703  & 389  \\ 
&ilu    				& 325  & 609  & 566  & 534  & 294  \\ 
&mg     				& 5    & 5    & 5    & 5    & 5    \\ 
\hline 
\multirow{5}{*}{$512 \times 512$}
&none   				& 3271 & 17685 & 18777 & 16890 & 6172  \\ 
&jacobi 				& 3429 & 17721 & 18777 & 16873 & 6105  \\ 
&sor    				& 1043 & 2699  & 2987  & 2510  & 1120  \\ 
&ilu    				& 764  & 2050  & 2045  & 1596  & 828   \\ 
&mg     				& 5    & 5     & 5     & 5     & 5     \\ 
\hline 
\end{tabular} 
\label{tab:pc}
\end{center}
\end{table}

In addition, it is interesting to look at the effect of the time step $\tau$.
The unpreconditioned GMRES method  was used with the grid of $256 \times 256$ unknowns.
From Table \ref{tab:tau} we see that the number of iterations decreases with $\tau$.
\begin{table}[!h]
\begin{center}
\caption{The number of iterations for various $\tau$}
\begin{tabular}{ccccccc}
\hline 
\multirow{2}{*}{$\eta$} &\multirow{2}{*}{$\tau$} & \multicolumn{5}{c}{$\xi$}\\
 						&  	   &$-10$ & $-1$ & $0.0$ & $1$ & $10$  \\ \hline
\multirow{5}{*}{0.01} 
							   &0.01 & 864  & 866  & 866  & 866  & 868 \\
							   &0.1  & 3329 & 3324 & 3323 & 3321 & 3306 \\
                               &1    & 4710 & 4686 & 4682 & 4678 & 4629 \\
                               &10   & 4913 & 4887 & 4882 & 4877 & 4822 \\
                               &100  & 4936 & 4907 & 4903 & 4899 & 4844 \\
                               \hline
\multirow{5}{*}{0.1} 
							   &0.01 & 892  & 924  & 926  & 928  & 933  \\
							   &0.1  & 3203 & 3420 & 3412 & 3396 & 3015 \\
                               &1    & 4413 & 4726 & 4689 & 4652 & 3828 \\
                               &10   & 4568 & 4914 & 4884 & 4830 & 3919 \\
                               &100  & 4687 & 4934 & 4903 & 4850 & 3939 \\
                               \hline
\multirow{5}{*}{1} 
							   &0.01 & 958  & 1413 & 1444 & 1446 & 1085 \\
							   &0.1  & 1313 & 3927 & 3953 & 3661 & 1614 \\
                               &1    & 1445 & 4790 & 4790 & 4289 & 1689 \\
                               &10   & 1455 & 4896 & 4893 & 4368 & 1695 \\
                               &100  & 1456 & 4908 & 4904 & 4374 & 1696 \\
                               \hline
\end{tabular} 

\label{tab:tau}
\end{center}
\end{table}

\section{Parallel implementation}\label{s6}
The parallel formulation is based on the domain decomposition methods. 
The main idea of these methods is to divide the original computational domain 
into subdomains.  A separate processor, which is identified by its rank, is assigned to each subdomain
in order to perform the computations. 
For inter processor communications the Message Passing Interface (MPI) is used. 

The systems of linear equations are solved by the parallel version of 
the preconditioned GMRES algorithm.
In our computations the none, bjacobi (doing the ILU-factorization of a local part of the matrix at each processor)  
and multigrid  preconditioners were used.  
The calculations were performed on the grid $ 512 \times 512 $ at $ \eta = 1 $.

The parallel code was run on a cluster of North--Eastern Federal University. 
The cluster consists of four computing nodes, each node has 
two quad-core processors Intel Xeon E5450 (3.00 GHz)  with 16 Gb RAM.

The results of the parallelization efficiency  of computations
are given in Table \ref{tab:par}. 
The table shows the estimation of computational costs, 
since the number of iterations is almost independent of the number of running processes.

\begin{table}[!h]
\begin{center}
\caption{Computation time in seconds, np --- number of processes, pc --- preconditioner}
\begin{tabular}{ccccccc}
\hline \multirow{2}{*}{np}  & \multirow{2}{*}{pc} & \multicolumn{5}{c}{$\xi$}\\ 
 &  & -10 & -1 & 0  & 1  & 10  \\ \hline
\multirow{5}{*}{16} 
							   &none    & 41.26  & 226.68  & 241.78  & 202.67  & 72.21  \\
                               &bjacobi & 11.48  & 31.16   & 31.69   & 24.94   & 12.66  \\
                               &mg      & 3.06   & 3.17    & 3.10    & 3.10    & 3.14   \\
                               \hline
\multirow{5}{*}{8} 
							   &none    & 57.35  & 304.48  & 236.72  &  294.56 & 108.07  \\
                               &bjacobi & 13.55  & 43.54   & 28.92   &  3502   & 12.46  \\
                               &mg      & 2.38   & 2.39    & 3.64    &  2.85   & 2.91  \\
                               \hline
\multirow{5}{*}{4} 
							   &none    & 86.28  & 593.00  & 583.87  & 538.89  & 220.66  \\
                               &bjacobi & 29.08  & 86.85   & 88.96   & 69.80   & 34.30  \\
                               &mg      & 4.99   & 4.15    & 3.38    & 4.84    & 4.92  \\
                               \hline
\multirow{5}{*}{2} 
							   &none    & 216.75  & 1328.398 & 1361.5  & 883.42  & 471.60 \\
                               &bjacobi & 75.55   & 177.12   & 145.78  & 146.34  & 79.90  \\
                               &mg      & 5.9     & 5.87     &  7.98   &  5.90   & 7.95   \\
                               \hline 
\multirow{5}{*}{1} 
							   &none    & 315.07  & 1686.16 & 1799.68 & 1612.7  & 590.24 \\
                               &bjacobi & 73.33   & 197.30  & 216.68  & 200.49  & 86.07  \\
                               &mg      & 9.23    & 8.55    & 8.57    &  8.56   & 8.54   \\
                               \hline                                                                             
\end{tabular} 
\label{tab:par}
\end{center}
\end{table}

\section{Conclusions}\label{s7}
\begin{enumerate}
 \item  The basic features of the pressure problem  
 associated with the non-selfadjoint operator are considered for multiphase flows in porous media.
 \item It was found that the computational cost of solving the model pressure problem does not depend 
 strongly on $ \xi $, more pronounced dependence is on the physical parameter $ \eta $. 
 This means that the number of iterations depends basically on the various properties 
 of the phases than on the value of external sources.
 \item Parallel computations have been performed using standard techniques with various preconditioners.
\end{enumerate}

\end{document}